\documentclass[12pt,twoside]{article}

\usepackage{latexsym}
\usepackage{epsfig} \epsfverbosetrue
\usepackage{headerfo}
\usepackage{amsfonts}

\newcounter{theorem}
\newcommand{\newsection}[1]{{\setcounter{theorem}{0} \section{#1}}}
\newtheorem{Theorem}{Theorem}[section]
\newtheorem{Proposition}[Theorem]{Proposition}
\newtheorem{Lemma}[Theorem]{Lemma}
\newtheorem{Corollary}[Theorem]{Corollary}

\newcommand{\Z}{\mathbb Z}
\newcommand{\R}{\mathbb R}
\newcommand{\C}{\mathbb C}

\newif\ifpdf
\ifx\pdfoutput\undefined
  \pdffalse
\else
  \pdfoutput=1
  \pdftrue
\fi

\usepackage{latexsym}
\ifpdf
  \usepackage[pdftex]{graphicx}
\else
  \usepackage{graphicx}
  
\fi

\chardef\aa=64

\begin{document}

\renewcommand{\author}{L.P.~Bos and A.~Brudnyi\\  
Department of Mathematics and Statistics\\
University of Calgary \\
Calgary, Alberta\\
Canada T2N 1N4 \\
\medskip
and \\
\medskip
N.~Levenberg\\
Department of Mathematics \\
University of Auckland \\
Private Bag 92019, Auckland, New Zealand \\}

\newcommand{\stitle}{Markov Inequalities on Exponential Curves}
 
\renewcommand{\title}{A Tangential Markov Inequality on
Exponential Curves}

\renewcommand{\date}{\today}
\flushbottom
\setcounter{page}{1}
\pageheaderlinetrue
\oddpageheader{}{\stitle}{\thepage}
\evenpageheader{\thepage}{\stitle}{}
\thispagestyle{empty}
\vskip1cm
\begin{center}
\LARGE{\bf \title}
\\[0.7cm]
\large{\author}
\end{center}
\vspace{0.3cm}
\begin{center}
\date
\end{center}
\begin{abstract}
We show that on the curves $y=e^{t(x)}$ where
$t(x)$ is a fixed polynomial, there holds a tangential Markov inequality
of exponent four. Specifically, for the real interval $[a,b]$
there is a constant $C$ such that
\[ \max_{x\in [a,b]}\left|{d\over dx}P(x,e^{t(x)})\right|\le
C(\hbox{deg}(P))^4\max_{x\in [a,b]}|P(x,e^{t(x)})| \]
for all bivariate polynomials $P(x,y).$
\end{abstract}

\setcounter{section}{0}

\newsection{Introduction} Recently [BLMT1, BLMT2, Ba1, Ba2, 
Br1, Br2, FN, RY] there has been considerable
interest on extending the classical Markov/Bernstein inequalities bounding
the derivatives of polynomials to higher dimensional cases. To be precise,

\noindent {\bf Definition}. Suppose that $M$ is a smooth ($C^1$)
manifold
in $\R^d$ (with or without boundary). 
We say that $M$ admits a tangential Markov 
inequality
of exponent $\ell$ if 
there is a constant $C>0$ such that for all
polynomials $P\in\R[x_1,\cdots,x_d]$ and points $a\in M$
$$|D_TP(a)|\le C(\hbox{deg}(P))^\ell ||P||_M.$$

Here $D_TP$ denotes any (unit) tangential derivative of $P$ and
$||P||_M$ is the supremum norm of $P$ on $M.$\par
\medskip

One remarkable fact is that such tangential Markov inequalities 
of exponent {\it one} characterize $M$ being algebraic.

\begin{Theorem} ([BLMT1]). Suppose that $M$ is a smooth
($C^\infty$) compact $m-$dimensional submanifold of $\R^d.$  
If $M$ admits a tangential
Markov inequality of exponent $\ell<1+1/m,$ then $M$ is algebraic
(i.e. a subset of an algebraic variety of the same dimension). Conversely,
if $M$ is algebraic, then it admits a tangential Markov inequality of
exponent $\ell=1.$
\end{Theorem}

It is natural then to investigate the existence of a tangential
Markov inequality of exponent $\ell>1$ and its connection
to the algebraicity of $M.$
For instance, in [BLMT2] it is shown that the
non-singular part of the
{\it singular} algebraic curve segment
$(x,x^r),$ $0\le x\le1,$ with $r=q/p>1$ in lowest terms, admits
a tangential Markov inequality of exponent $\ell=2p$ and this
is best possible. However, not all analytic manifolds
admit a tangential Markov inequality.
In [BLMT1] there is given an example of an {\it analytic}
curve in $\R^2$ which does not admit a tangential Markov
inequality of {\it any} exponent $\ell.$ Specifically, the
curve segment $(x,f(x)),$ $0\le x\le1$ with $f(x)$ given by a
gap series $f(x)=\sum_{k=0}^\infty c_k x^{2^{k!}}$ with $c_k\ge0$ and
being convergent on $[0,1],$ has this property. 

Thus the main question is to characterize manifolds which admit a tangential
Markov inequality. Up till now it was not even clear whether 
or not the existence of a Markov
inequality of {\it some} exponent meant that $M$ was (possibly singular)
algebraic. In this short note we show that this is not the case by showing that
the analytic curve segments $(x,e^{t(x)}),$ $a\le x\le b$ do admit a tangential
Markov inequality of exponent $\ell=4.$ Precisely, we prove

\begin{Theorem}
Suppose that $t(x)$ is a fixed algebraic polynomial. Then for every
interval $[a,b]$ there exists a constant $C=C(a,b)$ such that
\[ \max_{x\in [a,b]}\left|{d\over dx}P(x,e^{t(x)})\right|\le
C(\hbox{deg}(P))^4\max_{x\in [a,b]}|P(x,e^{t(x)})| \]
for all bivariate polynomials $P(x,y).$
\end{Theorem}

\newsection{Siciak Type Extremal Functions}

Suppose that $K\subset \C^d$ is compact and let
\[ {\cal P}_K=\{P\in\C[z_1,\cdots,z_d]\,:\, \|P\|_K\le1\,\,
\hbox{and}\,\,\hbox{deg}(P)\ge1\}.\]
(Here $\|P\|_K$ denotes the uniform norm of $P$ on $K.$)
The, by now classical,
Siciak extremal function ([S], but see also the excellent monograph [K]),
is defined, for $z\in \C^d,$ as 
\[\Phi_K(z)=\hbox{sup}\{|P(z)|^{1/{\rm deg}(P)}\,:\,P\in {\cal P}_K\}.\]
This function has proven to be a highly useful tool in the theory
of analytic functions of several variables, and is one of the main
ingredients in the proof of Theorem 1 above. However, the polynomials
restricted to the curve $w=e^{t(z)}$ are not themselves polynomials and hence
we need to generalize appropriately.

Suppose then that $\{V_n\}_{n\ge1}$ is a collection of increasing
subsets of $C(\C^d),$ i.e., $V_n\subset V_{n+1}\subset 
C(\C^d),$ $n=1,2,\cdots,$
and let $\displaystyle{V=\bigcup_{n=1}^\infty}V_n.$ For $f\in V$ we
set
\[ {\rm deg}(f)=\inf\{n\,:\,f\in V_n\}.\]
Further, let
\[\phi\,:\,\Z_+\to \R\]
be a given function. We define a generalized Siciak  extremal function
for the family $V$ (replacing polynomials) with exponent function $\phi$
(replacing ${\rm deg}(P)$) to be
\[ \Phi_K(z;V,\phi)=\sup_{f\in V\atop \|f\|_K\le1} |f(z)|^
{1/\phi({\rm deg}(f))}\qquad (z\in \C^d).\]
If there is no risk of confusion then we will drop the dependence on $V$ 
and $\phi.$

In our situation we take
\begin{equation}
V_n=\{f(z)=P(z,e^{t(z)})\,:\,P\in \C[z_1,z_2]\,\,\hbox{and}\,\,
{\rm deg}(P)\le n\}\subset C(\C^2)
\end{equation}
and
\begin{equation} 
\phi(n)=n^2.
\end{equation}

\noindent {\bf Remark.}
We emphasize that in the classical case $\phi(n)=n^1,$ but here
we must, as we shall see, use the higher power $n^2.$ This is a
crucial point.

\medskip Now, in the course of the proof of Lemma 1, page 120, of
[Bak], there is given an estimate on exponential polynomials
(i.e. for $t(z)=z$) due to Tijdeman[T]
which we may restate in the following form.
For fixed $f\in V$ 
and $z_0\in \C$ let $\displaystyle{M(R)=\max_{|z-z_0|\le R}|f(z)|}.$

\begin{Lemma}  Suppose that $t(z)=z.$  Then, 
there are constants $c_1,c_2>0$ such
that for all $f\in V$ and $R_2\ge R_1>0,$
\[ M(R_2)\le M(R_1)e^{c_1R_2n}\left({R_2\over R_1}\right)^{c_2\phi(n)} \]
where $n={\rm deg}(f)$ and $\phi$ is given by (2).
\end{Lemma}
\noindent {\bf Proof.} See [Bak], \S2 of Chapter 12. $\Box$

\noindent {\bf Remark.} Inequalities of the form of Lemma 2.1 are classically
referred to as Bernstein-Walsh inequalities. In case $R_2=2R_1$ they
are also sometimes called {\it doubling} inequalities.\par

\medskip Taking $1/n^2$ powers of this
estimate gives an immediate estimate for the Siciak extremal
function for $K$ a disk.

\begin{Corollary} Suppose that $t(z)=z$ and
that $K=\{z\in \C\,:\,|z-z_0|\le R_1\}$
(with $R_1>0$).
Then, there is a constant 
$\displaystyle{C=C(R_1,R_2)\le e^{c_1R_2}{R_2\over R_1}}$ such that
\begin{equation}
\Phi_K(z;V,\phi)\le C
\end{equation}
for all $|z-z_0|\le R_2.$
\end{Corollary}

Baker proves Lemma 2.1 in order to exploit the classical Jensen
inequality to obtain 
a bound on the valency of exponential polynomials (see [Bak, T] for
details). 
Further, Roytvarf and Yomdin[RY] prove a converse to this 
for general analytic functions which they
make use of in their
theory of Bernstein classes.

Perhaps the most succinct statement of this converse is
given by Brudnyi[Br2, Lemma 3.1] which we paraphrase as

\begin{Lemma} Suppose that $f$ is analytic on
$K=\{z\in \C\,:\,|z-z_0|\le R_2\}$ and assumes
no value there more than $p$ times.  Then,
for $0<R_1\le R_2$
there is a constant $C=C(R_1,R_2)$ such that
\[M(R_2)\le C^pM(R_1).\]
\end{Lemma}

Thus to establish that the extremal function is bounded for
polynomials restricted to the curve $w=e^{t(z)}$ we need
only give an appropriate bound on the valency of such functions,
which we do next.

\begin{Theorem}  Suppose that $t(z)$ is a fixed polynomial
of degree $d,$
and that $K=\{z\in \C\,:\,|z-z_0|\le R\}$ with $R\ge R_0>0.$  
There is a constant $C=C(t,K,R_0)$ such that
the valency in $K$ of every $f\in V_n$
is bounded, 
\[ p\le Cn^2.\]
More precisely, we have
\[ p\le A(dn^2+dn\left({R\over R_0}\right)^d
\sup_{|z-z_0|\le R_0}|t(z)|)\]
where $d=\hbox{deg}(t)$ and $A\le4$ is an absolute constant.
\end{Theorem}

\noindent {\bf Proof.} For simplicity we take $z_0=0.$ The idea
of the proof is to reduce to the case of $t(z)=z.$

So suppose that $P(z,w)$ is a bivariate polynomial of degree at
most $n.$  Then
\[P(z,e^{t(z)})=\sum_{k=0}^n a_k(z)e^{kt(z)}\]
for some polynomials $a_k(z)$ with ${\rm deg}(a_k)\le n-k.$
We would like to bound the number of zeros in $K$ of
$P_c(z,e^{t(z)}):=P(z,e^{t(z)})-c.$
For $w\in \C$ let $z_1,z_2,\cdots,z_d$ denote the $d$ zeros of $t(z)-w$ 
(repeated according to multiplicity) and set
\[ {\tilde P}_c(w)=\prod_{j=1}^d P_c(z_j,e^{t(z_j)}).\]
We claim that, in fact,
\[{\tilde P}_c(w)=\sum_{k=0}^{nd}b_k(w)e^{kw}\]
for some polynomials $b_k(w)$ with ${\rm deg}(b_k)\le n.$ To
see this, first note that $t(z_j)=w$ for each $j$ and
hence 
\[{\tilde P}_c(w)=\prod_{j=1}^d\left\{\sum_{k=0}^n
{\tilde a}_k(z_j)e^{kw}\right\}\]
where ${\tilde a}_0(z)=a_0(z)-c,$ ${\tilde a}_k(z)=a_k(z),$ $k\ge1.$
Upon expanding we see that
\[{\tilde P}_c(w)=\sum_{k=0}^{nd}b_k(z_1,\cdots,z_d)e^{kw}\]
where the $b_k$ are {\it symmetric} polynomials in the roots
$z_1,\cdots,z_d$ of degree at most $nd.$ But then the
$b_k$ must be polynomials in the 
elementary symmetric functions $\sigma_i(z_1,\cdots,z_d),$
i.e., in the coefficients of $t(z)-w.$
These are all constant, except for the constant coefficient
which is a function of $w.$ Further, the elementary symmetric
function $\sigma_d(z_1,\cdots,z_d)=\prod_{j=1}^dz_j$ is already
of degree $d$ and hence the degree of $b_k,$ as a polynomial
in $w,$ must be at most $nd/d=n.$

To continue, for $R\ge R_0$ set
\[ R_1:=\left({R\over R_0}\right)^d\sup_{|z|\le R_0}|t(z)|\]
so that, by the classical Bernstein inequality,
\[ \sup_{|z|\le R}|t(z)|\le R_1.\]
In other words, $|z|\le R$ $\Longrightarrow$ $|t(z)|\le R_1.$
Then, for $c\in\C,$
\begin{eqnarray*}
&&\#\{z\,:\, |z|\le R\,\,\hbox{and}\,\,P(z,e^{t(z)})=c\} \cr
&&\le\#\{z\,:\,|t(z)|\le R_1\,\hbox{and}\,\,P(z,e^{t(z)})=c\} \cr
&&\le\#\{w\,:\,|w|\le R_1\,\,\hbox{and}\,\,{\tilde P}_c(w)=0\} \cr
&&\le A(dn^2+dnR_1) \cr
\end{eqnarray*}
by Lemma 1 of [Bak, p. 120]. The result follows. $\Box$

In summary, Corollary 2.2 also holds for general polynomials
$t(z).$

\begin{Corollary}
Suppose that $K=\{z\in \C\,:\,|z-z_0|\le R_1\}$
(with $R_1>0$).
Then there is a constant 
$\C=C(R_1,R_2)$ such that
\begin{equation}
\Phi_K(z;V,\phi)\le C
\end{equation}
for all $|z-z_0|\le R_2.$ 
\end{Corollary}

\medskip\noindent {\bf Remark.} In the classical case the Siciak
extremal function has been much studied in the context of plurisubharmonic
functions. We do not know to what extent this carries over to
our generalized Siciak extremal functions, but for our purposes here
we need only know that it is bounded. \par

\medskip Corollary 2.5 gives a bound for the extremal function
of a complex disk, but since we wish to prove a Markov inequality
for a {\it real} segment of the curve $y=e^{t(x)},$ we will need to know
that the Siciak extremal function is bounded for a real segment.  Suppose
then that $x_0\in\R$ and let $B_R(x_0)$ be the {\it complex} disk
of radius $R$ centred at $x_0$ and let $I_R(x_0)$ denote the 
{\it real} segment, centred also at $x_0$ with radius $R.$

\begin{Theorem} Suppose that $V$ is a set of entire functions and is
such that 
\[ \Phi_{B_{R_1}(x_0)}(z;V,\phi)\le C(R_1,R_2)\quad \forall |z-x_0|\le R_2,
\,\,R_1\le R_2.\]
Then there is an associated constant $C'=C'(R_1,R_2)$ such that
\[ \Phi_{I_{R_1}(x_0)}(z;V,\phi)\le C'(R_1,R_2)\quad \forall |z-x_0|\le R_2,
\,\,R_1\le R_2.\]
\end{Theorem}

\noindent {\bf Remark.} Note, this holds for {\it general} $V$ and $\phi.$ 
Of course, later we will apply it to generalized exponential polynomials and
$\phi(n)=n^2.$ Moreover, we will show that $\phi(n)=n^2$ is best possible.

\medskip\noindent {\bf Proof.} The proof is an ingenious application,
due to the second author A.~Brudnyi[Br3], of a version of
a generalization of a classical Lemma of Cartan, due to Levin[L],
which we begin by quoting.

\medskip\noindent{\bf Theorem.} ([L, p. 21])
{\it Suppose that $f(z)$ is entire and such that $f(0)=1.$ Let
$0<\eta<3e/2$ and $R>0.$ Then,
inside the disk $|z|\le R,$ but outside a family of disks the
sum of whose radii is not greater than $4\eta R,$ we have
\[\ln\,|f(z)| >-H(\eta)\ln\, M(2eR) \]
for 
\[ H(\eta)=2+\ln\, {3e\over 2\eta}. \] }

\medskip Now consider $g\in V_n$ and let $z^*\in \C,$ $|z^*-x_0|=R_1,$ be
such that $\displaystyle{|g(z^*)|=\max_{|z-x_0|\le R_1}|g(z)|.}$

We apply Cartan's Lemma with $z^*$ the origin to
\[ f(z)={g(z)\over g(z^*)}, \]
$R=2R_1$ and $\eta=1/16.$  Then, there is a set of disks $\{D_i\}$
of radii $r_i$ such that $\sum_i r_i\le 4{1\over16}R=R_1/2<R_1$
and that, inside the disk $B_{2R_1}(z^*)\supset B_{R_1}(x_0),$
but {\it outside} the disks $D_i,$
\begin{equation}
\ln\,|f(z)|>-(2+\ln(24e))\ln\, M_f(4eR_1).
\end{equation}
But, since the sum of their radii is strictly less than $R_1,$
the disks $D_i$ cannot possibly cover all of $I_{R_1}(x_0),$ i.e.,
there is a point $x\in I_{R_1}(x_0)$ such that (5) holds.  For this
point then, by the defintion of $z^*,$
\[\ln\left({|g(x)|\over \max_{|z-x_0|\le R_1}|g(z)|}\right)
>-(2+\ln(24e))\ln\left( {\max_{|z-z^*|\le 4eR_1}|g(z)|\over
\max_{|z-x_0|\le R_1}|g(z)|}\right) \]
and hence
\[\ln\left({\max_{|x-x_0|\le R_1}|g(x)|\over \max_{|z-x_0|\le R_1}|g(z)|}\right)
>-(2+\ln(24e))\ln\left( {\max_{|z-z^*|\le 4eR_1}|g(z)|\over
\max_{|z-x_0|\le R_1}|g(z)|}\right). \]
Reversing the sign in the equality, we get
\begin{eqnarray*}
\displaystyle{\ln\left({\max_{|z-x_0|\le R_1}|g(z)|
\over\max_{|x-x_0|\le R_1}|g(x)| }\right)}&<&
\displaystyle{(2+\ln(24e))\ln\left( {\max_{|z-z^*|\le 4eR_1}|g(z)|\over
\max_{|z-x_0|\le R_1}|g(z)|}\right)} \\
&\le&\displaystyle{(2+\ln(24e))\ln\left( {\max_{|z-x_0|\le (4e+1)R_1}|g(z)|\over
\max_{|z-x_0|\le R_1}|g(z)|}\right)}.
\end{eqnarray*}
Consequently,
\begin{eqnarray*}
\displaystyle{{\max_{|z-x_0|\le R_1}|g(z)|
\over\max_{|x-x_0|\le R_1}|g(x)|}} &\le&
\displaystyle{\left({\max_{|z-x_0|\le (4e+1)R_1}|g(z)|\over
\max_{|z-x_0|\le R_1}|g(z)|}\right)^{2+\ln(24e)}} \\
&\le& \displaystyle{C(R_1,(4e+1)R_1)^{(2+\ln(24e))\phi(n)}}
\end{eqnarray*}
by the boundedness of the extremal function for complex disks. 
The result follows. $\Box$

\medskip We will use the next result to show that, for the curve $y=e^x,$ the
exponent function $\phi(n)=n^2$ is best possible. First some
notation. For $f(x),$ analytic in a neighbourhood of the origin,
let ${\rm ord}_0(f)$ denote the order of vanishing of
$f$ at $x=0,$ i.e., so that
\[ f(x)=cx^{{\rm ord}_0(f)}+\,\,{\rm higher\,\, order\,\, terms}\]
for some $c\neq0,$ near $x=0.$ For $W$ a vector space of functions, analytic in
a neighbourhood of the origin, we set
\[{\rm maxord}_0(W)=\sup_{f\in V} {\rm ord}_0(f),\]
the maximum order of vanishing for functions in $W.$ (We take
${\rm ord}_0(0)=0.$)

\begin{Proposition} Let $W$ be a vector space of functions,
analytic in a neighbourhood of the origin, such that
${\rm maxord}_0(W)<\infty.$  Set
\[ M_W(r)=\sup_{f\in W} {\max_{0\le x\le 2r}|f(x)|
\over \max_{0\le x\le r}|f(x)|}.\]
Then,
\[ \sup_{r\ge0} M_W(r)\ge 2^{{\rm maxord}_0(W)}.\]
\end{Proposition}
\noindent {\bf Proof.} For a fixed $0\neq f\in W,$
\[\lim_{r\to0}{\max_{0\le x\le 2r}|f(x)|
\over \max_{0\le x\le r}|f(x)|} =2^{{\rm ord}_0(f)},\]
as is easy to see.  Hence
\[ \sup_{r\ge0} M_W(r)\ge \sup_{f\in W} 2^{{\rm ord}_0(f)} =
2^{{\rm maxord}_0(W)}.\]
$\Box$

\medskip\noindent
Now, since $y=e^x$ is not algebraic, $V_n=\{P(x,e^x)\,:\,{\rm deg}(P)\le n\},$
is a vector space of dimension $\displaystyle{N:={n+2\choose 2}}.$

\begin{Lemma} 
\[ {\rm maxord}_0(V_n)=N-1.\]
\end{Lemma}

\noindent {\bf Proof.}
First note that, setting $f_{ij}(x)=x^ie^{jx},$ $i+j\le n$ positive
integers,
this is equivalent to the existence of constants $\alpha_{ij}$ such that
\[f^{(s)}(0)=0,\qquad 0\le s\le N-2\]
and
\[f^{(s)}(0)=1,\qquad s=N-1,\]
for $\displaystyle{f(x)=\sum_{ij}\alpha_{ij}f_{ij}(x).}$
But this in turn is equivalent to solving, for the $\alpha_{ij},$ the
linear system
\[\sum_{ij}\alpha_{ij}f_{ij}^{(s)}(0)=0\qquad 0\le s\le N-2,\]
\[ \sum_{ij}\alpha_{ij}f_{ij}^{(s)}(0)=1,\qquad s=N-1.\]
This can be done as the corresponding $(N-1)\times (N-1)$
determinant ${\rm det}[f_{ij}^{(s)}(0)]$ is the
Wronskian for the functions $\{f_{ij}=x^ie^{jx}\}_{i+j\le n}$ which
is non-zero due to the fact that these functions form a full set 
of linearly independent solutions of a certain constant coefficient
linear differential equation. $\Box$

\medskip Since $N=O(n^2),$ it follows from Proposition 2.7 and
Lemma 2.8 that, indeed, $\phi(n)=n^2$ is best possible in order 
for the extremal function $\Phi_{[0,r]}(x;V,\phi)$ to be bounded on
the interval $[0,2r].$

\medskip We now prove that a bounded extremal function implies the
existence of a Markov type inequality.  

\begin{Theorem} Suppose that $[a,b]\subset\R$ and that
the extremal function
$\Phi_{[a,b]}(z;V,\phi)$ is 
bounded inside the ellipse ${\cal E}_{a,b,c}\subset\C$
having foci at $a$ and $b$ and major axis $[a-c,b+c]$ for 
some $c>1/(\phi(1))^2.$ Then
there is a Markov inequality, in the sense that there exists a
constant $C>0$ such that
\[ \|f'(x)\|_{[a,b]}\le C(\phi(n))^2\|f\|_{[a,b]} \]
for all $f\in V_n.$
\end{Theorem}
\noindent {\bf Proof.} We make use of the so-called {\it relative}
extremal function of pluripotential theory. $C$ will denote
a generic constant; it need not have the same value in all
of its instances. 

Suppose that
$\Omega\subset\C^d$ is an open set and that $E\subset\Omega$ is
compact. Then the relative extremal function is defined as
\begin{equation}
u_{E,\Omega}(z):=\sup\{v(z)\,:\,v\in{\cal PSH}(\Omega),\,\,
v\le 1\,\,{\rm on}\,\, \Omega,\,\,v\le 0\,\,{\rm on}\,\,E\}.
\end{equation}
Here ${\cal PSH}(\Omega)$ denotes the plurisubharmonic functions on
$\Omega$ (see [K] for the precise defintion of this class of
functions). It is known that for $\Omega={\cal E}_{a,b,c}\subset \C$
and $E=[a,b]\subset\Omega,$
\begin{equation}
u_{E,\Omega}(z) ={\ln|w/r+\sqrt{(w/r)^2-1}|\over
\ln|R/r+\sqrt{(R/r)^2-1}|}
\end{equation}
where $w:=z-(a+b)/2,$ $r=(b-a)/2$ is the ``radius'' of the interval
$[a,b]$ and $R=r+c$ is that of $[a-c,b+c].$

Now, for non-constant $f\in V,$ with the same meanings 
of $\Omega$ and $E$ as above, let
\[ g(z)={ \ln|f(z)|-\ln\|f\|_E\over \ln\|f\|_\Omega-\ln\|f\|_E}.\]
Then $g\in {\cal PSH}(\Omega)$ and, clearly, $g\le1$ on $\Omega$ and
$g\le0$ on $E.$ Hence $g$ is a competitor in the definition of
the relative extremal function (6) and so
$g(z)\le u_{E,\Omega}(z).$ It follows that
\[ { \ln|f(z)|-\ln\|f\|_E\over \ln\|f\|_\Omega-
\ln\|f\|_E}\le u_{E,\Omega}(z) \]
and that
\begin{eqnarray}
\displaystyle{\ln\left({|f(z)|\over \|f\|_E}\right)}&\le&
\displaystyle{ \ln\left({\|f\|_\Omega\over\|f\|_E}
\right)u_{E,\Omega}(z)} \nonumber\\
&\le&\displaystyle{\ln\left(\sup_{z\in\Omega}\Phi_E(z;V,\phi)^{\phi({\rm deg}(f))}\right)
{\ln|w/r+\sqrt{(w/r)^2-1}|\over
\ln|R/r+\sqrt{(R/r)^2-1}|} }\nonumber\\
&\le&C\phi({\rm deg}(f))\ln|w/r+\sqrt{(w/r)^2-1}|,\,\,z\in\Omega
\end{eqnarray}
by the boundedness of the extremal function.

We now bound the derivative of $f\in V$ by means of the
Cauchy Integral Formula. Suppose then that $x\in [a,b]$ and
let $\Gamma_s\subset\C$ be the circle of radius $s$ centred
at $x.$ Then
\[f'(x)={1\over2\pi i}\int_{\Gamma_s} {f(z)\over (z-x)^2}\,dz\]
so that
\begin{equation}
|f'(x)|\le {1\over s}\max_{z\in \Gamma_s}|f(z)|.
\end{equation}
But then by (8), for $s$ sufficiently small so that
$\Gamma_s\subset\Omega,$
\[|f'(x)|\le {1\over s} \max_{z\in \Gamma_s}
|w/r+\sqrt{(w/r)^2-1}|^{C\phi({\rm deg}(f))} \|f\|_E\]
where $w$ and $r$ have the same meanings as previously stated.
Now, elementary calculations (cf. Lemma 1.1 of [BLMT1]) reveal that
there is a constant $C$ such that for $z-x\le s,$
\[ |w/r+\sqrt{(w/r)^2-1}|\le 1+C\sqrt{s}.\]
Hence,  
\[|f'(x)|\le {1\over s}
(1+C\sqrt{s})^{C\phi({\rm deg}(f))} \|f\|_E\]
and the result follows by taking 
\[ s={1\over (\phi({\rm deg}(f)))^2}.\]
$\Box$

\bigskip
Putting together Corollary 2.5, Theorem 2.6 and Theorem 2.9,
we obtain Theorem 1.2.

\newsection{Concluding Remarks}

Using methods similar to those presented in this paper, one
can show that $y=f(x)$ admits a tangential Markov inequality
of some exponent (which can be effectively estimated) if
$\displaystyle{f(x)=\sum_{i=1}^kp_i(x)e^{q_i(x)}}$
is a generalized exponential polynomial. Here $q_i,p_i
\in\R[x]$ are fixed polynomials. We conjecture that,
more generally, $y=f(x)$ admits a tangential Markov
inequality of some exponent if $f(x)$ satisfies
a linear differential equation 
\[y^{(d)}+a_1(x)y^{(d-1)}+\cdots+a_{d-1}(x)y'+a_d(x)y=0\]
with polynomial coefficients.

\bigskip

\centerline {\bf REFERENCES}
\vskip10pt
\begin{description}

\item{[Bak]} Baker, A., {\bf Transcendental Number Theory},
Cambridge U. Press, 1975.

\item{[Ba1]} Baran, M., {\it Bernstein theorems for
compact subsets of $\R^n$}, {\bf J. Approx. Theory}, 
69 (1992), 156--166.

\item{[Ba2]} Baran, M., {\it Bernstein theorems for
compact subsets of $\R^n$ revisited}, {\bf J. Approx. Theory}, 
79 (1994), 190-198.

\item{[BLMT1]} Bos, L., Levenberg, N., Milman P. and Taylor, B.A.,
{\it Tangential Markov Inequalities Characterize Algebraic
Submanifolds of $\R^n,$} {\bf Indiana J. Math.}, Vol. 44,
No. 1 (1995), pp. 115--138.

\item{[BLMT2]} Bos, L., Levenberg, N., Milman P. and Taylor, B.A.,
{\it Tangential Markov Inequalities on Real Algebraic Varieties},
{\bf Indiana J. Math.}, Vol. 47, No. 4 (1998), 1257--1271.

\item{[Br1]} Brudnyi, A., {\it Inequalities for entire functions},
in {\bf Entire Functions in Modern Analysis. 
Israel Mathematical Conferences  
Proceedings}, M. Sodin ed., to appear.

\item{[Br2]} Brudnyi, A., {\it Local inequalities for
plurisubharmonic functions}, {\bf Annals of Math.}, 149 (1999),
511--533.

\item {[Br3]} Brudnyi, A., {\it Bernstein type
inequalities for quasipolynomials}, to appear in
{\bf J. Approx. Theory}.

\item{[FN]} Fefferman, F. and Narasimhan, R., {\it A Local
Bernstein Inequality on Real Algebraic Varieties}, 
{\bf Mat. Z.}, 223 (1996), 673--692.

\item{[FY]} Francoise, J.-P. and Yomdin, Y., {\it Bernstein
Inequalities and Applications to Analytic Geometry and
Differential Equations},  {\bf J. Funct. Anal.} 146 (1997), 185--205.

\item{[K]} Klimek, K., {\bf Pluripotential Theory}, Oxford, 1991.

\item{[L]} Levin, B. Ja., {\bf Distribution of Zeros of Entire
Functions}, Trasnlations of Math. Monographs, Vol. 5, AMS, 1964.

\item{[RY]} Roytvarf, N. and Yomdin, Y., {\it Bernstein Classes},
{\bf Ann. Inst. Fourier}, Vol 47, No. 3 (1997), 825--858.

\item{[S]} Siciak, J., {\it On some extremal functions and their
applications in the theory of analytic functions of several
variables}, {\bf Trans. Amer. Math. Soc.}, 105 (1962), 322--357.

\item{[T]} Tijdeman, R., {\it On the number of zeros of general
exponential polynomials}, {\bf Indag. Math.}, Vol. 37 (1971),
1--7.
\end{description}

\end{document}